\documentclass{notices}
\usepackage{amsfonts,amssymb,amsmath,amscd}

\usepackage{latexsym}
\usepackage{graphics}
\usepackage{xspace}
\usepackage{amssymb}
\usepackage{psfrag}
\usepackage{epsfig}
\usepackage{amsthm}
\usepackage{pst-all}

\newcommand{\pr}{\mathbb{P}}

\newcommand{\R}{\mathbb{R}}

\newcommand{\G}{\mathbb{G}}

\newcommand{\mb}[1]{\mbox{\boldmath $#1$}}

\newcommand{\ignore}[1]{\relax}

\newtheorem{theorem}{Theorem}[section]

\newtheorem{Defi}[theorem]{Definition}

\newcommand{\ER}{Erd{\"o}s-R\'{e}nyi }
\newcommand{\Erdos}{Erd{\"o}s}

\newcommand{\Mezard}{M\'{e}zard}

\title{
Turing in the shadows of Nobel and Abel: an algorithmic story behind two recent prizes
}

\author{
 David Gamarnik\affil{
   David Gamarnik is a professor of operations research at the Massachusetts Institute of Technology. His email address is gamarnik@mit.edu.
   Funding from NSF grant CISE  2233897 is gratefully acknowledged.
    }
 }

\begin{document}

\maketitle

%\section{}

The 2021 Nobel Prize in physics was awarded to Giorgio Parisi  ``for the discovery of the interplay of disorder 
and fluctuations in physical systems from atomic to planetary scales,'' and the 2024 Abel Prize in mathematics was awarded 
to Michel Talagrand  ``for his groundbreaking contributions to probability theory and functional analysis, with outstanding 
applications in mathematical physics and statistics.'' What remains largely absent in the popular descriptions of these prizes, however,
is the profound contributions the works of both individuals have had to the 
field of \emph{algorithms and computation}. The ideas first developed by Parisi  and his collaborators 
relying on remarkably precise physics intuition,
and later confirmed by Talagrand and others by no less remarkable mathematical techniques, have revolutionized the way we think algorithmically
about optimization problems involving randomness. This is true  both in terms of  the existence of fast algorithms for some optimization problems, but also 
in terms of our persistent failures of finding such algorithms for some other optimization problems.
\ignore{
This important contribution stayed relatively in the shadows of the popular description of the prizes.
}

The goal of this article is to highlight these developments and
explain how the ideas pioneered by Parisi and Talagrand 
have led to a remarkably precise characterization of which optimization problems admit fast algorithms, versus those which do not,
and furthermore to explain why this characterization holds true.

\section{Optimization of random structures}  %The most embarrassing algorithmic open problem in probabilistic combinatorics}
The works of Parisi and Talagrand -- which were devoted to understanding a mysterious object in statistical physics
known as spin glass -- offered an entirely novel way of understanding algorithmic successes and failures 
in tackling optimization  problems involving randomness.  This progress was propelled by understanding 
the ``physics'' properties of the underlying problems, namely the geometry of the solution space. We will illustrate
these notions using three  examples, all of which  fit into  a general framework of optimization problems of  the form
\begin{align}
\eta_{\rm OPT} \triangleq \max_{\sigma\in\Theta} H(\sigma,\mathcal{R}). \label{eq:ground-state}
\end{align}
%\ignore{\textcolor{red}{[Isn't $\eta_{\rm OPT}$ a random variable of $\mathcal R$?]} }
Here $H$ is the function to be maximized, which we call the ``Hamiltonian'' for reasons to be
explained later,
$\Theta$ is the space of potential solutions, and $\mathcal{R}$ is a random variable describing the source of randomness.
Our three examples all exhibit a persistent Computation-to-Optimization gap, namely a gap between the optimal value $\eta_{\rm OPT}$
and the best guaranteed value that can be computed by known efficient algorithms -- where ``efficient'' means ``polynomial-time'' in the parlance of 
the theory of computational complexity.

\subsection*{Example I. Largest cliques and independent sets in random graphs}
The first of our three examples   is the problem of finding a largest clique in a dense \ER graph
and the related  problem of finding a largest independent set in a sparse \ER graph.

Imagine a club with $n$ members where every pair of members are friends or non-friends with equal likelihood
$1/2$, independently across all  ${n\choose 2}$ member pairs. What is the size of a largest group of members who are all friends with each other? This quantity is known as the 
clique number of the graph whose nodes are the members and whose edges correspond to the pairs of friends.  This club random friend pairings is just a vernacular description of a model of a random graph known as \ER graph,
denoted $\G(n,1/2)$. Paul \Erdos ~introduced this model in 1948 as a very effective tool to study and solve problems in the field of extremal combinatorics. The question of estimating 
the largest clique size of an \ER graph was solved successfully in the 1970s and has a remarkably precise answer: with high probability
(namely with probability approaching $1$ as $n\to\infty$, denoted as ``whp''), the size is $2(1+o(1))\log_2 n$. %~\cite{alon2004probabilistic}.  
Here $o(1)$ denotes a term which converges to zero
as $n\to\infty$. This is an instance of our optimization problem \eqref{eq:ground-state} in which $\Theta$ is the set of all subsets of the node
set $[n]:=\{1,\ldots,n\}$ of the graph, $\mathcal{R}$ is the random graph itself, and $H(\sigma,\mathcal{R})$ is the cardinality of $\sigma$
when $\sigma$ is a clique of $\mathcal{R}$, and $H=0$ otherwise.

The asymptotic value $2(1+o(1))\log_2 n$ has a very simple intuition. Fix an integer value $x$ and let  $Z(x)$ denote a random variable
corresponding to  the number of subsets of $[n]$  which have cardinality $x$ and are cliques. By the linearity of expectation, the expected value
of this random variable is ${n\choose x}\left({1\over 2}\right)^{x\choose 2}$. By straightforward computation this quantity
converges to zero when $x$ is asymptotically larger than $2\log_2 n$, which provides a simple proof of the upper bound.
The matching lower bound requires more work involving estimation of the variance of $Z(x)$ which is now standard material in many textbooks on 
random graphs.  \ignore{f.e. \cite{alon2004probabilistic}.}

What about the problem of efficiently finding  a largest clique algorithmically in at most polynomial time, 
when the graph is given to the algorithm designer as an input? 
Here and elsewhere we adhere to the gold standard of algorithmic efficiency when
requiring that the running time of the algorithm grows at most  polynomially in the size of the input. 
That is we require that the algorithm belongs to  $\mathcal{P}$ -- the class  of algorithm whose running time grows at most
 polynomially in the size of the problem description, which in our case means it grows at most polynomially in $n$. 
 Note that the largest clique problem can be trivially solved 
in time $n^{O(\log n)}$ by   checking all subsets with cardinality $2(1+o(1))\log_2 n$, but this growth rate is super-polynomial and thus
this exhaustive search algorithm is not in $\mathcal{P}$.

This is where our story begins. In 1976 Richard Karp -- one of the founding fathers of the modern algorithmic complexity theory -- presented a very simple algorithm that computes a clique whose size is guaranteed to be (roughly) half-optimal, namely $(1+o(1))\log_2 n$~\cite{karp1976probabilistic}. 
Given that his algorithm was not particularly ``clever,'' Karp challenged the community to produce an algorithm with a better performance guarantee,
with the stipulation that the algorithm still remains in  $\mathcal{P}$.   Remarkably,
in the  almost-50 years since the problem was introduced
there has been no successful improvement of Karp's rather na\"{i}ve algorithm, 
including notably the algorithm based on Markov Chain Monte Carlo updates, as shown by Jerrum in 1992. %~\cite{jerrum1992large}.
This persistent failure further motivated Jerrum to introduce 
one of the most widely studied problems in random structures and  high dimensional statistics, namely the ''infamous'' Hidden Clique Problem.
The problem of finding large cliques in $\G(n,1/2)$ thus  exhibits Computation-to-Optimization gap.

Notably, while  improving the half-optimality for the clique problem appears algorithmically hard, 
this hardness is only \emph{evidential}
but not  (at least not yet)  formal. Namely, 
we do not have any formal evidence of even conditional algorithmic hardness for this problem 
say by assuming the validity of the ``$\mathcal{P}\ne \mathcal{NP}$'' conjecture. This is in stark contrast  with the 
worst-case setup. Namely, consider the problem of finding a largest clique for \emph{arbitrary} graph, in particular for graphs which are
 not necessarily generated according
to the $\G(n,1/2)$ distribution. For this version of the  problem it is known that no polynomial time algorithm can produce a clique within 
 any multiplicative factor $n^{1-\epsilon}$ from optimality (in particular much larger than factor $1/2$)
 unless $\mathcal{P}\ne \mathcal{NP}$.
We see that, in particular, the half-optimality gap emerging in  the random case  does not align with the classical algorithmic complexity theory based 
on the $\mathcal{P}\ne \mathcal{NP}$ conjecture.

A very similar story takes place   when the underlying  model is a  sparse random graph  defined next. Fix a constant $d$. Consider a random graph
denoted by $\G(n,d/n)$ in which every pair of nodes is an edge 
with probability $d/n$, independently across all pairs. 
We now switch to the mirror image of the clique problem, known as the independent set problem
defined as follows. A subset of nodes $I\subset [n]$ is called an independent set 
if it contains no edges. Namely for all $i,j\in I$ the pair $(i,j)$ is not an edge. 
Frieze proved~\cite{FriezeIndependentSet} that the size of the largest independent set is $\eta_{d,\rm OPT} \triangleq 2(1+o_d(1)){\log d\over d}n$ 
whp as $n\to\infty$, where here $o_d(1)$ denotes a constant vanishing
in $d$ as $d$ increases.

Consider now the problem of identifying a largest independent set algorithmically in polynomial time. 
Using the formalism described in  \eqref{eq:ground-state}, $\Theta$ is the set of all subsets of the node
set $[n]:=\{1,\ldots,n\}$ of the graph, $\mathcal{R}$ is the random graph $\G(n,d/n)$  itself, and $H(\sigma,\mathcal{R})$ is the cardinality of $\sigma$
when $\sigma$ is an independent set  in the graph  $\G(n,d/n)$,   and $H=0$ otherwise.
 An algorithm similar to Karp's clique algorithm computes an independent set whose size is roughly half-optimal, namely
$(1+o_d(1)){\log d\over d}n$, also  whp. No improvement of this algorithm is known and no formal hardness is known either. 
Thus, this problem is also evidentially hard and exhibits Computation-to-Optimization gap.

\subsection*{Example II. Ground states of spin glasses}
Our second example is in fact a probabilistic model of spin glasses themselves. 
 Consider a tensor of random variables 
$J=\left(J_{i_1,\ldots,i_p}, 1\le i_1<i_2<\cdots< i_p\le n\right)\in \R^{n\otimes p}$,
where each entry $J_{i_1,\ldots,i_p}$ of the tensor has a standard normal distribution and all entries are probabilistically independent.
The optimization problem of interest is 
\begin{align}\label{eq:SpinGlassGroundState}
\eta_{n,p}\triangleq \max_{\sigma \in \{0,1\}^n} n^{-{p+1\over 2}}\langle J,\sigma^{\otimes p}\rangle,
\end{align}
where for every $\sigma\in \{-1,1\}^n$,
\begin{align}\label{eq:Jinnersigma}
 \langle J,\sigma^{\otimes p}\rangle=\sum_{1\le i_1<i_2<\cdots< i_p\le n}J_{i_1,\ldots,i_p}\sigma_{i_1}\sigma_{i_2}\cdots\sigma_{i_p}.
\end{align}
The optimal solution $\sigma^*$ of this problem (which is  unique up to a sign flip) is called 
the \emph{ground state} of the model using the physics jargon, and the optimal value is called the \emph{ground state value}. The normalization $n^{-{p+1\over 2}}$
is chosen for convenience (see below). In our formalism (\ref{eq:ground-state})
$\Theta=\{-1,1\}^n$, $\mathcal{R}$ is the tensor $J$, and $H(\sigma,\mathcal{R})$ 
is the associated
inner product (\ref{eq:Jinnersigma}).
The special case $p=2$ corresponds to the well-known and widely studied Sherrington-Kirkpatrick 
model. 

Notice how the signs of $J_{i_1,i_2}$ impact the optimal choices for $\sigma_i$.
When $J_{i_1,i_2}>0$, $\sigma_{i_1}$ and $\sigma_{i_2}$ prefer to have the same sign so that the product $J_{i_1,i_2}\sigma_{i_1}\sigma_{i_2}$
is positive. On the other hand, when $J_{i_1,i_2}<0$, $\sigma_{i_1}$ and $\sigma_{i_2}$ prefer to have 
opposite signs, so that the product $J_{i_1,i_2}\sigma_{i_1}\sigma_{i_2}$ is again positive. What is the best arrangements 
of signs to $\sigma_i$ which optimizes this tradeoff? This is the core
of the optimization problem, which turned out to be quite challenging. 

In 1980 Parisi solved the optimization problem (\ref{eq:SpinGlassGroundState}) by introducing a mathematically non-rigorous but remarkably powerful
physics style method called Replica Symmetry Breaking (RSB)~\cite{parisi1980sequence}. 
The importance of the RSB based technique is hard to overstate: it is the method of choice
for computing the ground state values and related quantities
for complex statistical physics models. Parisi's solution of the problem (\ref{eq:SpinGlassGroundState}) is as follows.
 Let $\mathcal{U}$ be the set of non-negative non-decreasing right-continuous functions $\mu$ on $[0,1]$ satisfying $\int_0^1\mu(t)dt<\infty$.
Consider the solution $\Psi_\mu(t,x), (t,x)\in [0,1)\times \R$ of the following PDE with boundary condition $\Psi(1,x)=|x|$:
\begin{align*}
\partial_t \Psi_\mu=-{p(p-1) \over 2} t^{p-2}\left(\partial^2_x\Psi_\mu+\mu\left(\partial_x \Psi\right)^2\right).
\end{align*}
Define $\mathcal{P}(\mu)=\Psi_\mu(0,0)-{p(p-1)\over 2}\int_0^1 t^{p-2}\mu(t)dt$.

\begin{theorem}
For every $p$  the ground state value satisfies
\begin{align}\label{eq:Parisi-limit}
\lim_n \eta_{p,n}= \inf_{\mu\in\mathcal{U}}\mathcal{P}(\mu)\triangleq \eta_{p, \rm OPT},
\end{align}
whp as $n\to\infty$.
\end{theorem}

The optimal solution $\mu_{p,\rm OPT}$ of the variational problem above is known to be unique and defines a 
cumulative distribution function of a probability
measure, which has been called since then the Parisi measure. The identity (\ref{eq:Parisi-limit}) was derived by Parisi heuristically using the RSB method.
It was turned into a theorem by 
Talagrand in 2006 in his breakthrough  paper~\cite{talagrand2006parisi}. 
(Thus one could say that the scientific value of the formula (\ref{eq:Parisi-limit}) is worthy of two major prizes!) Talagrand's approach relied
on a powerful interpolation bound developed earlier by Guerra~\cite{guerra2003broken}.

What about  algorithms? Here the goal is to find a solution $\sigma$ such that 
$n^{-{p+1\over 2}}\langle J, \sigma^{\otimes p}\rangle$ is as close to $\eta_{p, \rm OPT}$ as possible,
when the tensor $J$ is given as an input. Interestingly,  major progress on this question was achieved only very recently and it deserves a separate
discussion found in Section~\ref{section:Subag+Montanari}. Jumping ahead though, the summary of the state of the art algorithms is as follows.
For $p=2$ near optimum solutions can be found by a polynomial time
algorithm. On the other hand, when $p>2$  this is evidently not the case. In summary, the largest
value $\eta_{p,\rm ALG}$  achievable by known poly-time algorithms satisfies $\eta_{2,\rm ALG}=\eta_{2,\rm OPT}$, 
$\eta_{p,\rm ALG}<\eta_{p,\rm OPT}, p>2$.
Thus similarly to \textbf{Example I}, the model exhibits a  Computation-to-Optimization gap 
between the   optimal values and the values found by the best known polynomial time algorithms.

\subsection*{Example III. Satisfiability of a random K-SAT formula}
Our third and the final family of examples are randomly generated constraint satisfaction problems, specifically the so-called random K-SAT problem.
An instance of a K-SAT model is described by $n$ Boolean variables $x_1,\ldots,x_n$ (variables with $0/1$ values), and a conjunction
$\Phi\triangleq \wedge_{1\le j\le m} C_j$ of 
$m$ Boolean functions $C_1,\ldots,C_m$, each of which is a disjunction of precisely $K$ variables from the list $x_1,\ldots,x_n$ or their negation. 
Here is an example for $n=10$:  
$\Phi=\left( x_3\vee \neg x_7 \vee \neg x_8\right) \wedge \left( \neg x_1\vee \neg x_2 \vee  x_7\right) \wedge \left( \neg x_3\vee \neg x_7 \vee \neg x_9\right)$.
An assignment of variables $x_1,\ldots,x_n$ to values $0$ and $1$ 
 is called a satisfying assignment if $\Phi=1$ under this assignment. 
 For example, the assignment  $x_2=0,x_3=1,x_9=0$ (and any assignment for  the remaining $x_i$-s)
is a satisfying assignment. The formula $\Phi$ is called satisfiable if there exists at least one satisfying assignment and is called
non-satisfiable otherwise.

A random K-SAT formula is generated according to the following probabilistic construction. For each $j=1,\ldots,m$  independently across
$j$ we select $K$ variables $x_{j,1},\ldots,x_{j,K}$ uniformly at random from $x_1,\ldots,x_n$ and attach the negation sign $\neg$ to each of them 
with probability $1/2$, also independently
across $j=1,2,\ldots,m$ and $k=1,\ldots,K$. The resulting formula $\mb{\Phi}_{n,m}$ is a random formula whose probability distribution 
depends   $n$ and $m$ only.
In the context of the optimization problem (\ref{eq:ground-state}), 
$\Theta=\{0,1\}^n$, $\mathcal{R}$ is the description of the formula $\mb{\Phi}_{n,m}$, and $H(\sigma,\mathcal{R})$ 
is the number of clauses satisfiable by the assignment $\sigma\in\Theta$. In particular the formula is satisfiable if and only if  
$\max_{\sigma}H(\sigma,\mathcal{R})=m$.

Whether or not $\mb{\Phi}_{n,m}$ is  satisfiable or not is now an event, which we write symbolically as $\mb{\Phi}_{n,m}\in {\rm SAT}$. The probability of this event
$\pr\left(\mb{\Phi}_{n,m}\in {\rm SAT}\right)$ is our focus. There exists a (conjectural) critical threshold $c_{K,\rm SAT}$  such that
$\lim_n \pr\left(\mb{\Phi}_{n,cn}\in {\rm SAT}\right)=1$ when $c<c_{K,\rm SAT}$ and $=0$ when $c>c_{K,\rm SAT}$. Informally, $c_{K,\rm SAT}$
represents the largest fraction of clauses to variables for which the solution exists whp asymptotically in $n$. The existence of such a threshold
$c_{K,\rm SAT}$ was recently proven for large enough $K$~\cite{ding2022proof}. 

The constant  $c_{K,\rm SAT}$  has a very simple 
asymptotic expansion for large $K$: $c_{K,\rm SAT}=2^K\log 2(1+o_K(1))$, where $o_K(1)$ hides a function vanishing to zero as $K\to\infty$.
The intuition for this asymptotics is rather straightforward. Observe that any \emph{fixed} assignment of $x_i, 1\le i\le n$ is satisfiable with probability
$(1-1/2^K)^m$. Thus the expected number of satisfying solutions is $2^n(1-1/2^K)^m$. The value $2^K\log 2$ is precisely the critical value of $c$
at which this expectation drops from $\exp(\Theta(n))$ (exponentially large) to $\exp(-\Theta(n))$ (exponentially small). 

Turning to algorithms,  the best known one succeeds in finding
a satisfying assignment only when $c\le c_0 2^K {\log K \over K}\triangleq c_{K, \rm ALG}$~\cite{coja2010better}. 
Here $c_0$ denotes a universal  constant. Namely, 
 the best known algorithms succeed only  when the clause to variables density is nearly a factor $K/\log K$ (!) smaller 
than the  density for which the solutions are   guaranteed to exist whp. This is an even more dramatic  gap, than the half-optimality
gap seen earlier for cliques/independent sets in the  $\G(n,1/2)$ and $\G(n,d/n)$ models.
\ignore{
Quoting Dimitris Achlioptas, one
of the pioneers of the field in early 2000s: ''Solutions exist with  strongest probabilistic guarantees one could hope for; 
there are tons  of them, and
yet, we fail to find them. This is crazy!''.}
Random K-SAT is thus  another instantiation of the Computation-to-Optimization gap,
similar to the ones found in \textbf{Examples I, II}. 

In the remainder  of the article we will  explain 
the nature of these gaps. In particular, we will  explain the values  $1/2, \eta_{p, \rm ALG}$ and $c_{K, \rm SAT}$
by leveraging powerful ideas first developed towards understanding statistical physics of  materials called spin glasses.

\section{Physics matters. Spin glasses}
Statistical physics studies properties of matter using probabilistic  models which provide statistical  likelihoods one  finds physical systems in different states.
These likelihoods in their turn are linked with the energy of these states: the smaller  the energy, the bigger the likelihood.
One can apply this view point to our optimization problem $\max_\sigma H(\sigma,\cdot)$, and think of  $H(\sigma,\cdot)$ as  (negative) 
``energy'' (Hamiltonian) of the state $\sigma$. The ground state is then a  state with the lowest energy. 
The link between energies and likelihoods  is formalized by introducing the Gibbs measure
on the state space $\Theta$ which assigns probability measure $\langle \sigma \rangle =Z^{-1}\exp(\beta H(\sigma,\cdot))$ to each configuration
$\sigma$. Here $Z$ is the normalizing constant called the partition function. 
The parameter $\beta$ is called inverse temperature. At low temperature, namely large values 
of $\beta$, the Gibbs measure
becomes concentrated in near ground states: states with nearly lowest energy. The identity (\ref{eq:Parisi-limit}) was first derived for the case
of finite $\beta$ and then extended to the case $\beta=\infty$, corresponding to the ground state value.
The added value  of introducing  the Gibbs measure is that it allows to consider the entire space $\Theta$ 
of solutions in one shot where the solutions are stratified
by their energy levels. This is the so-called ``solution space geometry'' approach which turned out to be very fruitful for the analysis of algorithms.

Spin glasses are particularly complex examples of matters which exhibit properties of being liquid and solid at the same time. 
For a recent treatment of the subject see ''Spin glass theory and far beyond: Replica symmetry breaking after
 40 years'',  edited by Marinari, Parisi, Ricci-Tersenghi, Sicuro, Zamponi, Mezard; World Scientific 2023.
 Mathematically, 
spin glasses correspond to the fact that the associated Hamiltonian $H$ which is a sum $\sum_j H_j$ of many local Hamiltonians representing interactions
of nearby atoms, is ``frustrated''. Put differently, it is impossible to find a spin configuration $\sigma$ which simultaneously minimizes the energy
of each component $H_j$. This makes it very similar to problems we have considered in our examples. 
\ignore{Thus understanding the solution
space geometry of spin glasses promises to provide a valuable insight into the solution space geometry for our optimization problems.}

The ideas of repurposing the techniques from the theory of spin glasses  for the field of optimization of combinatorial structures goes back to
mid-eighties of the previous century. In particular, Fu and  Anderson in their  1986 
paper ``Application of statistical mechanics to NP-complete problems in combinatorial optimisation'' say the following  
``...we propose to study the general properties of such problems (combinatorial optimization) in the light of
... the structure of solution space, the existence and the nature of phase transitions...''.

\section{Clustering phase transition in random K-SAT}
The first fruitful and mathematically verified 
link between the solution space geometry and phase transition on the one hand, and the algorithmic tractability/hardness on the other hand
appeared in  concurrent works of Achlioptas and Ricci-Tersenghi~\cite{achlioptas2006solution} 
and \Mezard, Mora and Zecchina~\cite{mezard2005clustering} in the context
of the random  K-SAT model introduced earlier. Given a random formula $\mb{\Phi}_{n,m}$ let $\mathcal{SAT}_{n,m}\subset \{0,1\}^n$ denote the (random) set of
satisfying assignments. Replica Symmetry Breaking calculations, which were
 originally introduced to study  spin glasses,  suggest heuristically that 
the \emph{bulk} of this space should be highly clustered
at some large enough clause density $c<c_{K,\rm SAT}$. This is what these papers established in a mathematically rigorous  way. 
In particular, they proved the following theorem.

\begin{theorem}
For every $K\ge 8$ there exists $c_{K, \rm Clust}<c_{K,\rm SAT}$ such that
 the following holds when $c\in (c_{K, \rm Clust},c_{K,\rm SAT})$:  whp as $n\to\infty$ there exists a disoint partition 
 $\mathcal{SAT}_{n,m}=\left(\cup_{\ell} \mathcal{SAT}_\ell\right) \cup \mathcal{SAT}_o, \ell\ge 2,$ 
 such that the Hamming distance between any two distinct sets  $\mathcal{SAT}_{\ell_1}$
  $\mathcal{SAT}_{\ell_2}, \ell_1\ne \ell_2$ is  $\Omega(n)$, and $|\mathcal{SAT}_o |/ |\mathcal{SAT}_{n,m}|=\exp(-\Omega(n))$.
\end{theorem}

\begin{figure}
\begin{center}
\scalebox{.31}
%{\includegraphics{Images/Clustering-3.pdf}}
{\includegraphics{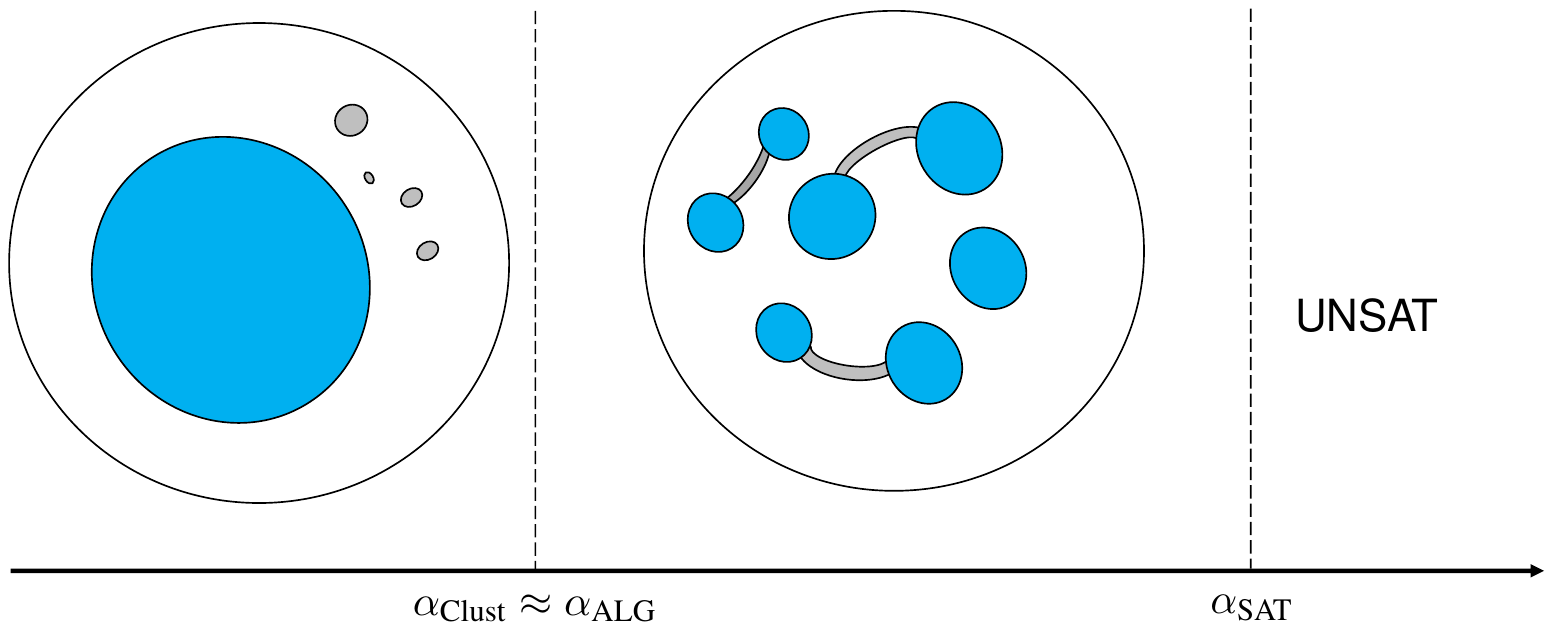}}
\end{center}
\caption{Clustering phase transition at $\alpha_{\rm Clust}\approx \alpha_{\rm ALG}$ . 
Below $\alpha_{\rm Clust}$ a bulk of the solution space (blue) is a giant connected set.
Above $\alpha_{\rm Clust}$ a bulk of the set is  partitioned into (blue) clusters, with exponentially small exceptions (grey)}
\label{fig:clustering}
\end{figure}

This is illustrated on Figure~\ref{fig:clustering}. The authors argued that the presence of such a 
 clustering picture suggests  algorithmic hardness of finding
satisfying assignments, in particular by means of some local search type algorithms. 
 Furthermore, as subsequent works showed, intriguingly,  the onset 
of this clustering  property phase transition coincides   with the evidential algorithmic hardness discussed earlier, asymptotically as $K$. More precisely,
it was verified that 
$c_{K, \rm Clust}=(1+o_K(1))2^K{\log K\over K}=(1+o_K(1))c_{K,\rm ALG}$! This was the first (though purely hypothetical) link between the solution 
space geometry  and the algorithmic  hardness. 

Very quickly a similar result was established also for the independent set problem in $\G(n,d/n)$ and several other similar problems: 
for all $c\in (1,2)$, independent
sets of size $c(1+o_d(1){\log d\over d}n$ (namely more than half optimum) are clustered, and are not clustered for all $c\in (0,1)$. 

What about the clustering property in spin glasses? It does take place and 
this is in fact what originally prompted Parisi to develop his Replica Symmetry Breaking approach.
Clustering  in spin glasses  is more intricate and comes in the form of nested hierarchy of sub-clusters arranged 
in a particular way  called ultrametric tree.
 We delay the discussion of this structure till 
Section~\ref{section:Subag+Montanari} when we discuss how this ultrametric structure of clusters guided successfully 
the creation of the state of the art algorithm for the problem (\ref{eq:SpinGlassGroundState}).

\section{Overlap Gap Property. A provable  barrier to algorithms}
As discussed earlier, the clustering property exhibited by random K-SAT and similar other problems  provided unfortunately only 
a hypothetical explanation of the algorithmic hardness, but not a proof of it. 
The property which turned out to be a provable barrier to algorithms was inspired by the  clustering picture,  but comes in the different
form of the Overlap Gap Property (OGP), which we now formally introduce in the context of the
problem (\ref{eq:ground-state}). For this purpose we assume that the solution space $\Theta$ is equipped with some metric $d$ (for example Hamming
distance when $\Theta$ is binary cube). 

\begin{Defi}\label{def:2-OGP}
An instance $H(\cdot,\mathcal{R})$ of the  optimization problem (\ref{eq:ground-state}) 
exhibits an OGP with parameters $\eta_{\rm OGP}\le  \eta_{\rm OPT}, \nu_1<\nu_2$ if for every
$\sigma,\tau\in \Theta$ satisfying $H(\sigma,\mathcal{R}),H(\tau,\mathcal{R})\ge \eta_{\rm OGP}$ it holds $d(\sigma,\tau)\notin (\nu_1,\nu_2)$.
\end{Defi}
This property has a  simple interpretation. The model exhibits OGP if the distance between 
every pair of ``near''-optimal solutions (solutions with values at least
$\eta_{\rm OGP}$) is either ``small'' (has a high overlap represented by $d(\sigma,\tau)\le \nu_1$) 
or is ``large'' (has a low overlap represented by $d(\sigma,\tau)\ge \nu_2$).  

We need to extend the definition above in two ways. First involving ensembles of instances $\mathcal{R}_1,\ldots,\mathcal{R}_m$, and second involving
overlaps of more than two solutions.

\begin{Defi}\label{def:e-m-OGP}
Fix a positive integer $M$ and a subset $\mb{\nu}\subset \Theta^M$. A family of problem instances $\mathcal{R}_1,\ldots,\mathcal{R}_T$
exhibits an OGP with parameters $(\eta_{\rm OGP},\mb{\nu})$ if for every $1\le \ell_1\le \ell_2\le\cdots \le \ell_M\le T$
and every $\sigma_1,\ldots,\sigma_M \in \Theta$ satisfying 
$H(\sigma_1, \mathcal{R}_{\ell_1}), \ldots, H(\sigma_M, \mathcal{R}_{\ell_M})\ge \eta_{\rm OGP}$, it holds 
$\left(\sigma_1,\ldots,\sigma_M\right)\notin \mb{\nu}$.
\end{Defi}
Definition~\ref{def:2-OGP} is recovered as a special case when $T=1, M=2$ and 
\begin{align}
\Theta^2 \supset \mb{\nu}=\{(\sigma,\tau): d(\sigma,\tau)\in (\nu_1,\nu_2)\}. \label{eq:2-ogp}
\end{align} 
We stress that the randomness $\mathcal{R}$ appearing in $H(\sigma_j, \mathcal{R}_{\ell_j})$ is index-specific, representing the fact that $\sigma_j$
is near-optimal (achieves value at least $\eta_{\rm OGP}$) specifically for instance $\mathcal{R}_{\ell_j}$.
While the intuition behind Definition~\ref{def:2-OGP} is fairly direct, the one for \ref{def:e-m-OGP} is rather opaque  and requires some elaboration,
to  be done  shortly. 

It turns out that the OGP does hold in our models and is a provable obstruction to \emph{classes} of algorithms. 
The first of these fact (in the context of Definition~\ref{def:2-OGP}) was already implicitly 
established in~\cite{achlioptas2006solution,mezard2005clustering} as a method for proving
the clustering property for the random K-SAT model. The fact though that OGP is a provable obstruction to algorithms 
 was established first in~\cite{gamarnik2014limits} in the context
of independent sets in sparse \ER graphs searchable by \emph{local} algorithms.  
 Loosely speaking an algorithm is local if for each node $i\in [n]$ the decision to  include it into an independent set is conducted entirely 
  based on a constant size graph-theoretic neighborhood around $i$, simultaneously for all $i$,  ignoring the rest of the graph.  
The result in~\cite{gamarnik2014limits} says the following.

\begin{theorem}\label{theorem:OGP-IS}
For every $\epsilon>0$ and sufficiently large $d$ there exist $0< \nu_1<\nu_2<1$ such that 
whp as $n\to\infty$ the set of independent sets in $\G(n,d/n)$ exhibits an OGP (in the sense of Definition~\ref{def:2-OGP})
with parameters $\left(\eta_{\rm OGP}\triangleq  \left({1\over 2}+{1\over 2\sqrt{2}}+\epsilon\right)\eta_{d,\rm OPT}; \nu_1n, \nu_2n \right)$. As a result, 
no local algorithm (appropriately  defined) exist for constructing an independent set with size larger than $\eta_{\rm OGP}$. 
\end{theorem}
Namely, every two independent sets with cardinality which is multiplicative factor ${1\over 2}+{1\over 2\sqrt{2}}$ close to optimality
are  either at  a ``small'' distance from each other (at most $\nu_1 n$) or at a ``large'' distance from each other (at least $\nu_2 n$). 
Theorem~\ref{theorem:OGP-IS}
means that beating factor ${1\over 2}+{1\over 2\sqrt{2}}$ optimality threshold for this problem will require going beyond the class of local algorithms. 
Local algorithms do in fact achieve
half-optimality (which recall is the evidential barrier for this problem).

The gap between half-optimality and  ${1\over 2}+{1\over 2\sqrt{2}}$ was  eliminated in~\cite{rahman2017local}  for the same class of local algorithms,  and later 
for a larger class of algorithms, see Theorem~\ref{theorem:IS-K-SAT} below. This was achieved  by turning to multi-overlaps
of solutions and resorting to  Definition~\ref{def:e-m-OGP}. A method to address broader classes of algorithms requires interpolating between
instances, that is having $T\ge 2$ in Definition~\ref{def:e-m-OGP},  as was introduced in~\cite{chen2019suboptimality}. 

The broadest currently known class of algorithms ruled out by OGP is the class of stable (noise insensitive) algorithms. We will not give one formal
definition of this class, as it comes in different context-dependent variants. Loosely speaking  the stability property 
means the following. We think of an algorithm as a 
map $\mathcal{ALG}: \mathcal{R}\to \Theta$ which maps instances $\mathcal{R}$ of  random input (graphs, tensors, etc) to solutions.
We say that the algorithm $\mathcal{ALG}$ is stable if small perturbation of the input $\mathcal{R}$ results in a small perturbation of the solution
$\mathcal{ALG}(\mathcal{R})$. Namely, if $\mathcal{R}\approx \mathcal{R}'$ then $\mathcal{ALG}(\mathcal{R})\approx \mathcal{ALG}(\mathcal{R}')$.
A local algorithm is an example of a stable algorithm: deleting or adding one edge of a graph will affect the decisions  associated only  with nodes which
contain this edge in their constant size neighborhood. This is so because the algorithm performs these decisions  concurrently (in parallel).
Thus  deleting or adding one edge changes the 
independent set produced by an algorithm in at most $O(1)$ many nodes.
Other classes of  algorithms proven to be successful for these types of problems, including 
Approximate Message Passing (see Section~\ref{section:Subag+Montanari}), stochastic localizations, low degree polynomials,
Boolean circuits with small depth, turn out to be stable as well. 

It was established in a series of followup works that OGP is a barrier to stable algorithms~\cite{gamarnik2024hardness},
\cite{gamarnik2022disordered}. For our problems outlined as \textbf{Examples I, III} it provides tight algorithmic bounds as described
in the following theorem, stated informally. The sequence of interpolated instances involving these examples is constructed as follows.
For our \textbf{Example I}, involving sparse  random graph $\G(n,d/n)$, consider any  ordering $\mathcal{P}$ 
of pairs $(i,j), 1\le i<j\le n$ and consider the process of resampling
each edge according to the same Bernoulli distribution with values $1,0$ with probabilities $d/n,1-d/n$, respectively, one pair at a time according to the
order $\mathcal{P}$. At each step $t$ we obtain a random graph $\G_t$ which is distributionally identical to $\G(n,d/n)$. After $T={n\choose 2}$
steps we obtain a fresh copy of of $\G(n,d/n)$ which is probabilistically independent from the starting copy. This is our interpolated sequence
$\mathcal{R}_1,\ldots,\mathcal{R}_T$ in the context of \textbf{Example I}. It is constructed similarly for the case of random K-SAT (\textbf{Example III}).

\begin{theorem}\label{theorem:IS-K-SAT}
For every large enough $d$ there exist a large enough $M$ and a subset $\mb{\nu}$
such  that whp the set of independent sets in the interpolation sequence 
$\G_t, 0\le t\le T={n\choose 2}$ of graphs $\G(n,d/n)$ exhibits an OGP with parameters $\eta_{\rm OGP}=1/2(1+o_d(1))$ and  $\mb{\nu}$. 
As a result, a largest independent set which can be constructed by stable algorithms is asymptotically $\eta_{\rm OGP}=1/2(1+o_d(1))$.

Similarly, for every large enough $K$ there exists a large enough  $M$ and a subset $\mb{\nu}$ 
such that whp the set of satisfying assignments  in the interpolation sequence
 $\mb{\Phi}_t$ of a random K-SAT formula  $\mb{\Phi}_{n,\alpha n}$ exhibits an OGP with parameters $\eta_{\rm OGP}=\eta_{\rm OPT}$
and  $\mb{\nu}$,
when $\alpha>\alpha_{K,\rm ALG}$. As a result the largest constraint-to-variables density for which stable algorithms
can find a solution of a random K-SAT problem is $\alpha_{K,\rm ALG}$. 

In particular, in both cases OGP provides a tight characterization of the power of stable algorithms.
\end{theorem}
For the second part of the theorem we set $\eta_{\rm OGP}=\eta_{\rm OPT}$ since for  constraint satisfaction problems
there is no notion of approximation to optimality, as we consider exclusively
satisfying assignments. Variants of the theorem above were proven 
in~\cite{rahman2017local, bresler2022algorithmic}.

\begin{figure}
\begin{center}
\scalebox{.31}{\includegraphics{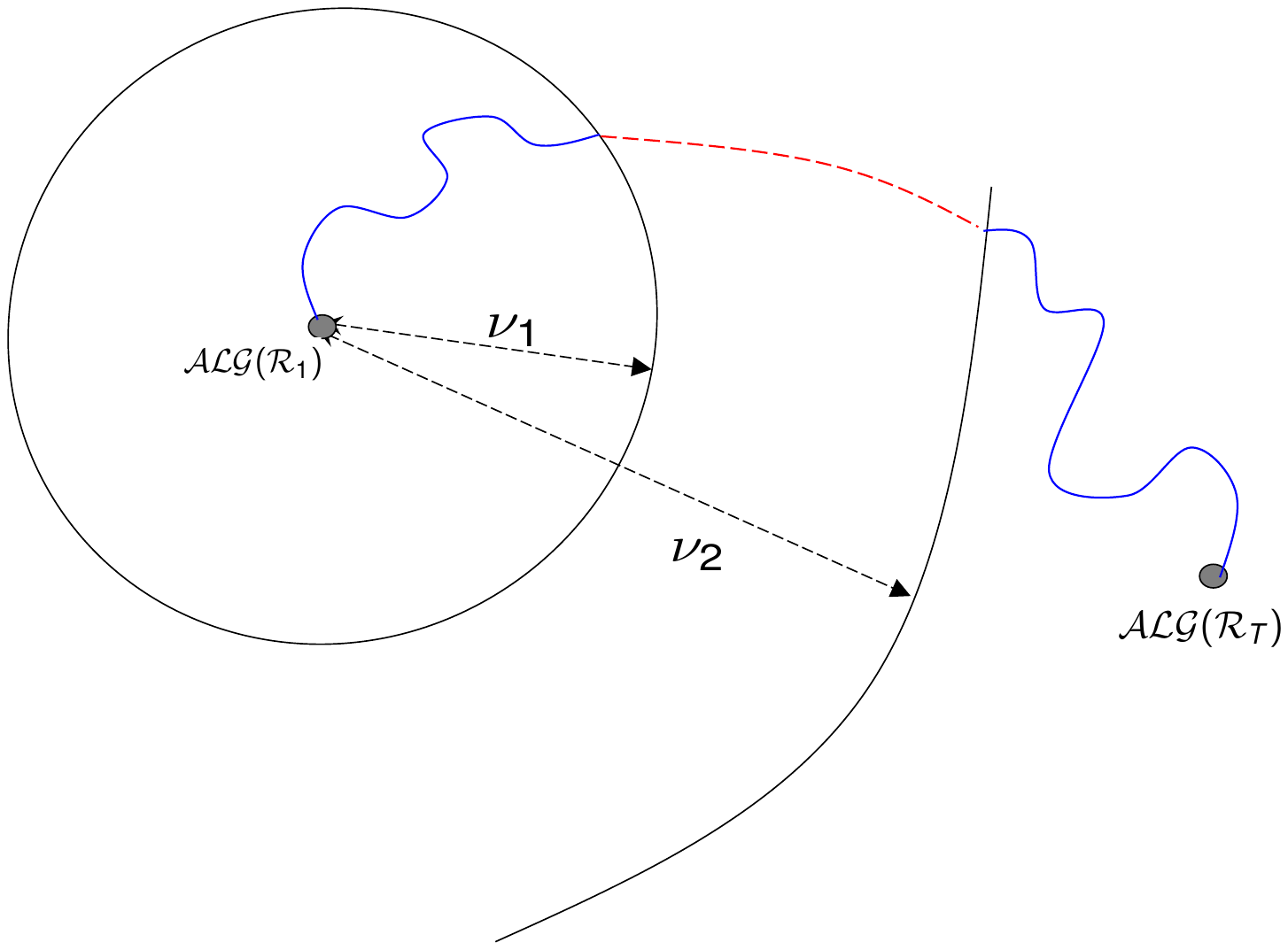}}
\end{center}
\caption{Overlap Gap Property is an obstruction to stable algorithms.
The algorithmic path (blue curve) has to jump (red section) from the small distance ($\le \nu_1$) region to the large distance 
($\ge\nu_2$) region in the space of solutions $\Theta$.}
\label{fig:OGP-stable}
\end{figure}

Why is OGP a barrier to stable algorithms? The intuitive  reason for this is a rather simple continuity type argument and best illustrated when $K=2$
in Definition~\ref{def:e-m-OGP}, and $\mb{\nu}$ is of the form (\ref{eq:2-ogp}). 
Suppose that in addition to the properties outlined in
Definition~\ref{def:e-m-OGP}, two additional properties hold. First, the sequence of instances $\mathcal{R}_1,\ldots,\mathcal{R}_T$ is ``continuous'',
in the sense that $\mathcal{R}_t\approx \mathcal{R}_{t+1}$ for all $t$. Second suppose 
that when $\ell_1=1$ and $\ell_K=\ell_2=T$, the case $d(\sigma_1,\sigma_2)\le \nu_1$
does not occur. Namely, while per Definition~\ref{def:e-m-OGP}, for every pair $\sigma_{t_1},\sigma_{t_2}$ of  solutions  with values $\ge \eta_{\rm OPT}$
in the interpolated
sequence $\mathcal{R}_t, t=1,\ldots,T$, both $d(\sigma_{t_1},\sigma_{t_2})\le \nu_1$ (small distance) and $\ge \nu_2$ (large distance) 
situations are possible, when the two 
solutions correspond to the beginning (that is $t_1=1$) and the end (that is $t_2=T$) of the interpolation sequence, only the latter  situation 
$d(\sigma_{t_1},\sigma_{t_2})\ge \nu_2$ (large distance) is possible. This can indeed be established in many special cases using the same techniques
which verify OGP. Then the argument establishing OGP as an obstruction to stable algorithms goes informally as follows. 
Consider implementing a stable algorithm $\mathcal{ALG}$ on a sequence of interpolated instances $\mathcal{R}_t, 1\le t\le T$,
resulting in a series of solutions $\mathcal{ALG}(\mathcal{R}_t)$. Since the sequence $\mathcal{R}_t$ is ``continuous'' (changes slowly as $t$ increases),
and the algorithm is stable, then $d\left(\mathcal{ALG}(\mathcal{R}_1),\mathcal{ALG}(\mathcal{R}_t)\right)$ changes ``continuously'' as well, as $t$ changes
from its initial value $d\left(\mathcal{ALG}(\mathcal{R}_1),\mathcal{ALG}(\mathcal{R}_1)\right)=0\le\nu_1$ to its final value
$d\left(\mathcal{ALG}(\mathcal{R}_1),\mathcal{ALG}(\mathcal{R}_T)\right)$ which by a property above must be $\ge \nu_2$. For any  $t$ by the 
OGP we have $d\left(\mathcal{ALG}(\mathcal{R}_1),\mathcal{ALG}(\mathcal{R}_t)\right)$ is either $\le \nu_1$ or $\ge \nu_2$. Therefore,  there must exist
$\tau$ such that $d\left(\mathcal{ALG}(\mathcal{R}_1),\mathcal{ALG}(\mathcal{R}_\tau)\right)\le \nu_1$ 
and $d\left(\mathcal{ALG}(\mathcal{R}_1),\mathcal{ALG}(\mathcal{R}_{\tau+1})\right)\ge \nu_2$. By triangle inequality this means
$d\left(\mathcal{ALG}(\mathcal{R}_\tau),\mathcal{ALG}(\mathcal{R}_{\tau+1})\right)\ge \nu_2-\nu_1$. But this contradicts stability of $\mathcal{ALG}$.
The argument is illustrated on Figure~\ref{fig:OGP-stable}

Recall now the clustering property discussed in the earlier section. 
Could  one build an algorithmic barrier type argument directly from the clustering property? It appears the answer
is no as evidenced in a different model called Symmetric Binary Perceptron model~\cite{GamarnikKizildagPerkinsXu2022}. In this model,
which also involves constraints-to-variables parameter $\alpha$, algorithms find a satisfying solution up to some value $\alpha_{\rm ALG}>0$
and the model exhibits an OGP above $\alpha_{\rm OGP}\approx\alpha_{\rm ALG}>0$. Yet, at the same the  model exhibits clustering at \emph{every}
positive value $\alpha>0$. 
It thus appears that geometric structures which are barriers to algorithms are 
more involved than clustering and have to rely on more complex structures 
such as multioverlaps $\mb{\nu}$.
We will see this in an even more pronounced form  in Section~\ref{section:Branching-OGP} in the context of spin glasses.

\section{Ultrametricity guides algorithms for spin glasses}\label{section:Subag+Montanari}
We now turn to the state of the art algorithms for spin glasses, our \textbf{Example II}, specifically algorithms for solving the optimization problem
(\ref{eq:ground-state}). A dramatic progress on this problem and a fairly complete understanding
of its algorithmic tractability  was achieved only recently in the works of Subag~\cite{subag2021following}, in the context of a related
so-called spherical spin glass model, and in Montanari~\cite{montanari2021optimization} in the context of precisely problem
(\ref{eq:ground-state}) in the special case $p=2$, namely the Sherrington-Kirkpatrick model. Montanari has constructed a polynomial
time algorithm which achieves any desired proximity $\eta_{p,\rm OPT}-\epsilon$ to optimality whp, assuming an unproven though widely believed
conjecture that the cumulative distribution function associated with the  Parisi measure $\mu_{p,\rm OPT}$ solving the variational problem
(\ref{eq:Parisi-limit}) is strictly increasing within its support.
Namely, the support of the measure induced by $\mu_{p,\rm OPT}$ is a connected subset of $[0,1]$. We will connect this property with OGP shortly. 

Soon afterwards, the result was extended to general $p$ in~\cite{el2021optimization} 
and a new algorithmic threshold $\eta_{p, \rm ALG}\le \eta_{p,\rm OPT}$ corresponding to the best known stable algorithms was identified.
This is summarized  in the  theorem below, stated informally. 
\ignore{
One  detail omitted in the statement of the theorem is the nature
of stability. Here the stability is understood in Lipschitz continuity sense.
Known algorithms achieving
the algorithmic thresholds, specifically the so-called Incremental Approximate Message Passing (IAMP) algorithm,  do satisfy this notion of stability.
}

\begin{theorem}[Informal]\label{theorem:eta-ALG-below}
Let $\mathcal{L}\supset \mathcal{U}$ be the set  of functions on $[0,1]$ with  total variation appropriately bounded.  Let 
$\eta_{p, \rm ALG}\triangleq \inf_{\mu\in\mathcal{L}}\mathcal{P}(\mu)\le \eta_{p, \rm OPT}=\inf_{\mu\in\mathcal{U}}\mathcal{P}(\mu)$. 
There exists
a polynomial time stable algorithm which 
constructs a spin configuration $\sigma_{\rm ALG}$ achieving $\approx \eta_{p, \rm ALG}$ whp as $n\to\infty$.
\end{theorem}

When $p=2$, it is conjectured that the variational problem $\inf_{\mu\in\mathcal{L}}\mathcal{P}(\mu)$ is solved
by a strictly increasing function, and thus modulo this conjecture  Montanari's algorithm achieves near optimality. 
It is known however~\cite{chen2019suboptimality} that this is not the case when $p\ge 4$ and $p$ is even. In this case 
\begin{align}\label{eq:Parisi-OGP}
\eta_{p, \rm ALG}=\inf_{\mu\in\mathcal{L}}\mathcal{P}(\mu)< \eta_{p, \rm OPT}=\inf_{\mu\in\mathcal{U}}\mathcal{P}(\mu).
\end{align}

\begin{figure}
\begin{center}
\scalebox{.33}{\includegraphics{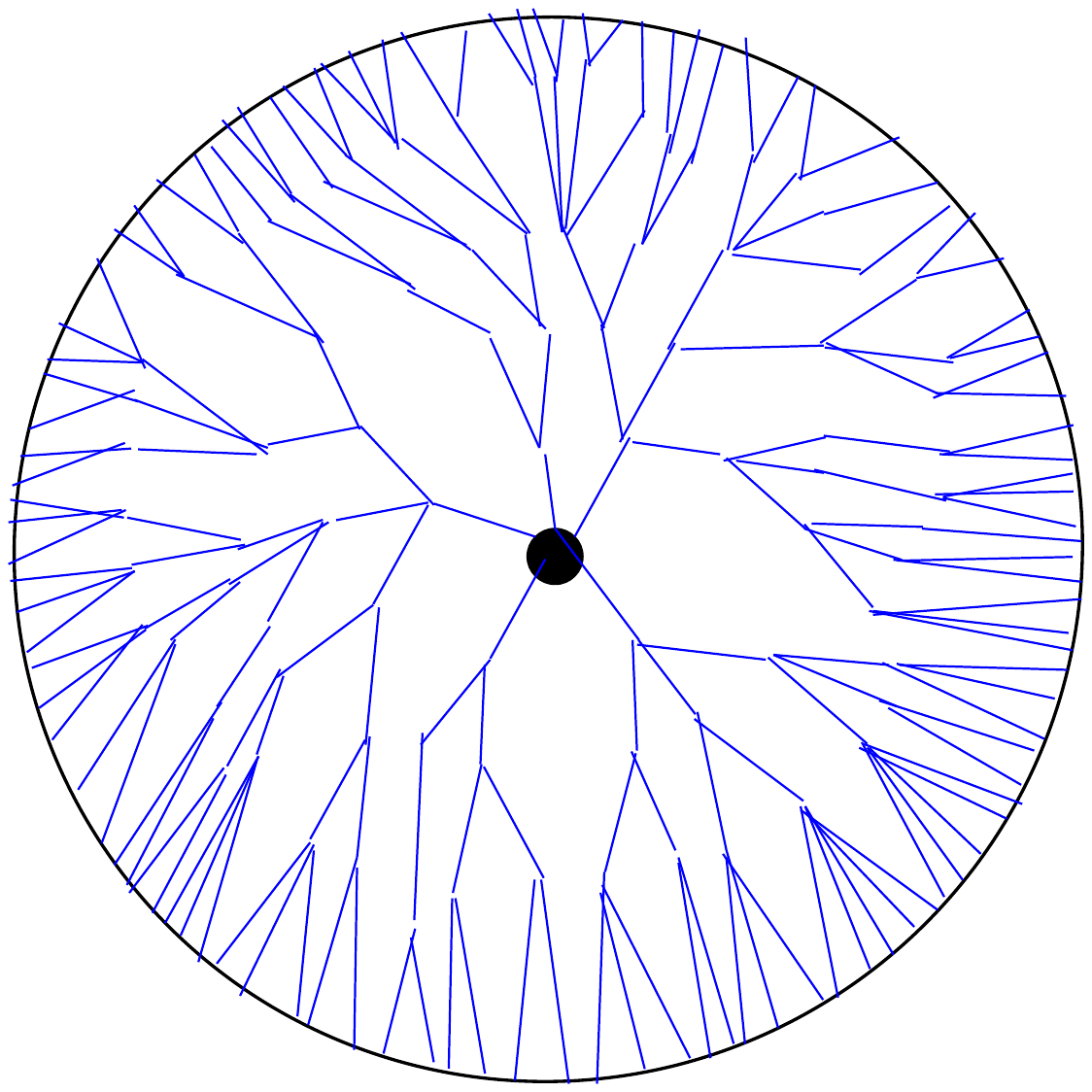}}
\end{center}
\caption{Ultrametric tree of solutions in $[-1,1]^n$.  End points of the curve correspond to corners of the cube. 
The Parisi measure has no gap in the support, reflected in the tree branching continously.}
\label{fig:UTree-no-gap}
\end{figure}

We now explain the algorithmic threshold $\eta_{p, \rm ALG}$ and the successful algorithm for achieving it. The Parisi
measure $\mu_{p,\rm OPT}$ is conjecturally a measure describing random overlaps of nearly optimal solutions in the following sense.
Fix a positive inverse temperature parameter $\beta$  and consider sampling two spin configurations $\sigma_1,\sigma_2$
from the associated Gibbs measure independently. Then $\mu_{p,\rm OPT}$ is the distribution of the inner product (overlap)
 $n^{-1}\langle \sigma_1,\sigma_2\rangle$ in the limit  $\beta\to\infty$.
Furthermore, suppose instead $m$ samples $\sigma_j, j\in [m]$ are generated according to the same Gibbs distribution. 
Then the distribution of the set of overlaps $n^{-1}\langle \sigma_{j_1},\sigma_{j_2}\rangle$ is believed to be described by distributions over set
of $m$  values  in $[0,1]$ arranged on a tree. As such it  satisfies a strong form of metric triangle inequality (with max  replacing the  addition operator) 
called ultrametricity, depicted on Figure~\ref{fig:UTree-no-gap}.  When the Parisi measure $\mu_{p,\rm OPT}$ is strictly increasing and there is no gap in its
support, the branching points on tree are located in a continuous way without gaps. It turns out that this coincides with the case 
$\inf_{\mu\in\mathcal{L}}\mathcal{P}(\mu)=\inf_{\mu\in\mathcal{U}}\mathcal{P}(\mu)$ implying
$\eta_{p, \rm ALG}=\eta_{p, \rm OPT}$, and this is the case when the algorithm reaches near optimality. The  algorithm
is constructed as a walk in the interior of $[-1,1]^n$ which starts from the center of the cube and  follows the ultrametric tree as a guide. Specifically, at every
stage it makes a small incremental step in the direction (a) maximally improving the objective function and (b) orthogonal to the previous direction.
Absence of gaps implies that the algorithm proceeds without interruption until hitting a corner of $[-1,1]^n$ provably reaching 
near optimality. In a sense, the ultrametric tree ``guides'' the successive steps of the algorithm.

When the Parisi measure is not strictly increasing, the authors show that the algorithm reaches the largest energy level achievable by ultrametric 
trees without gaps. This was shown to be the value of the extended variational problem 
$\inf_{\mu\in\mathcal{L}}\mathcal{P}(\mu)$ which per (\ref{eq:Parisi-OGP}) is strictly smaller than $\inf_{\mu\in\mathcal{U}}\mathcal{P}(\mu)$.

\section{Branching-OGP. Tight characterization of algorithmically solvable spin glasses}\label{section:Branching-OGP}
A final missing piece completing the algorithmic  tractability vs. hardness picture for spin glasses was delivered by Huang and Sellke~\cite{huang2022tight}.
When the Parisi measure $\mu_{p,\rm OPT}$ exhibits a gap in its support and as a result we have a strict inequality (\ref{eq:Parisi-OGP}), 
this inequality can be turned into an obstruction to algorithms along the lines of Theorems~\ref{theorem:OGP-IS} and~\ref{theorem:IS-K-SAT}.
Specifically, in the spirit of Definition~\ref{def:e-m-OGP}, let $\mb{\nu}_{p,\rm gap}$ be the set of $M$-tuples $(\sigma_j, j\in [M]) \in \Theta^M$
such that the set of the associated overlaps $\left(n^{-1}\langle \sigma_i,\sigma_j\rangle, i,j\in [M]\right)$ 
represents a continuous branching tree,  in some precise sense,   the kind we see in 
Figure~\ref{fig:UTree-no-gap}. Then the model exhibits an OGP above the value $\eta_{p,\rm ALG}$ with $\mb{\nu}$ serving as the multi-overlap 
obstruction. More specifically we have:

\begin{theorem}\label{theorem:b-OGP}
The optimization problem (\ref{eq:SpinGlassGroundState}) exhibits an OGP in the sense of Definition~\ref{def:e-m-OGP}
with parameters $\left(\eta_{p,\rm OGP}\triangleq \eta_{p,\rm ALG},\mb{\nu}_{p,\rm gap}\right)$.
This property is an obstruction to  stable algorithms in finding spin configurations with values larger than $\eta_{p,\rm ALG}$, whp as $n\to\infty$.

Hence (when combined with Theorem~\ref{theorem:eta-ALG-below}) $\eta_{p,\rm ALG}$ represents a tight algorithmic value achievable
by stable algorithms. 
\end{theorem}

As a corollary, a $p$-spin model can be solved to optimality by stable algorithms when $p=2$ (modulo the strict monotonicity conjecture discussed above).
On the other hand, when $p\ge 4$ is even,  the largest value achievable by stable algorithms is strictly sub-optimal. It is conjectured that this remains
the case for all  $p\ge 3$ regardless of the parity and has been verified for the case when  $p$ is sufficiently large. 
The multi-OGP obstruction set $\mb{\nu}_{p,\rm gap}$ 
is called Branching-OGP (B-OGP). It subsumes the obstruction sets $\mb{\nu}$ appearing in the context of 
Theorems~\ref{theorem:OGP-IS} and~\ref{theorem:IS-K-SAT} as special cases, and represents the most general known form of OGP-based obstructions. 

\begin{figure}
\begin{center}
\scalebox{.33}{\includegraphics{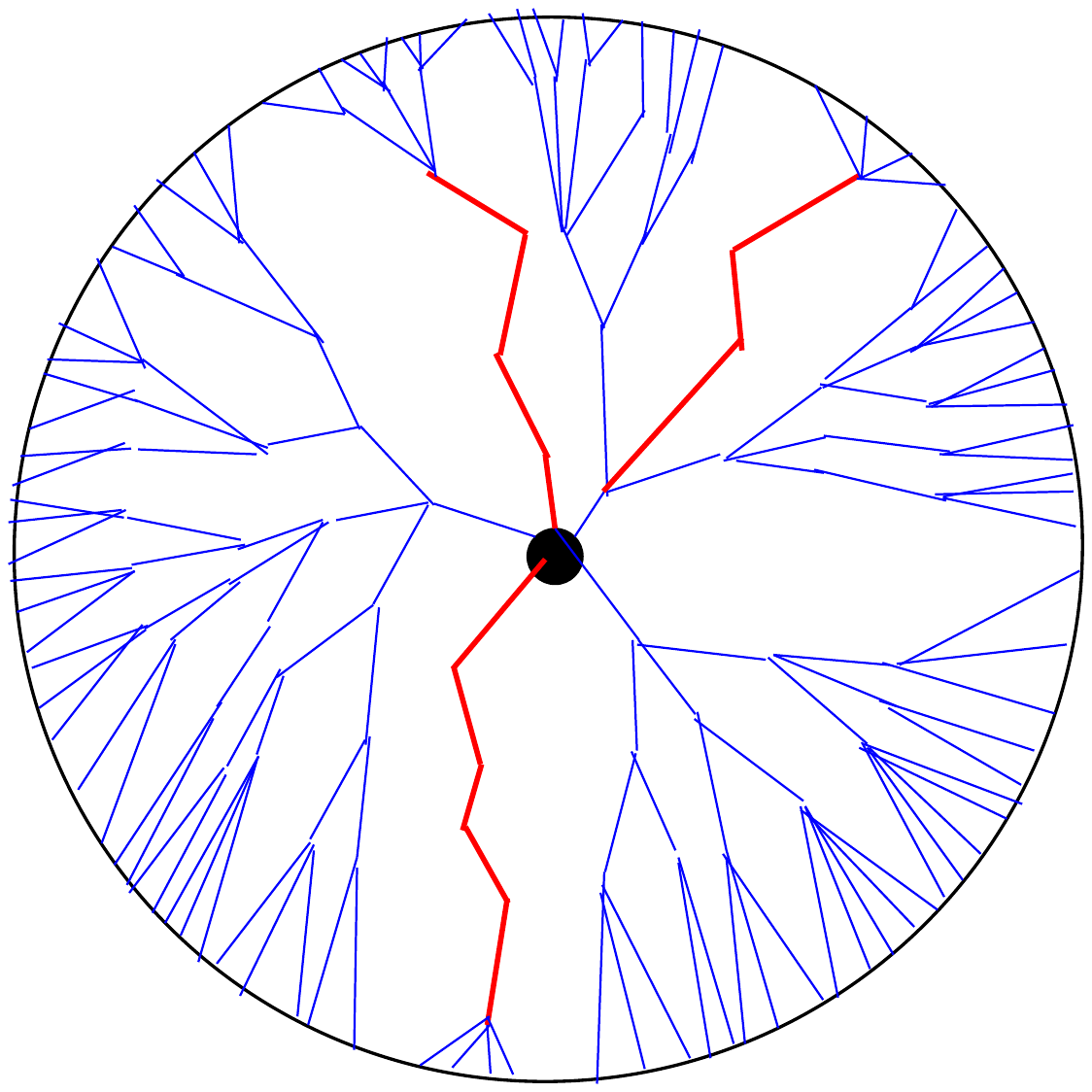}}
\end{center}
\caption{Ultrametric tree of solutions when the  Parisi measure has a gap in the support leading to gaps in branching locations.}
\label{fig:UTree-B-OGP}
\end{figure}

The proof of Theorem~\ref{theorem:b-OGP} follows a path similar to the ones used in Theorems~\ref{theorem:OGP-IS} and~\ref{theorem:IS-K-SAT} 
and is also based on the interpolation scheme between instances, using in a crucial way the  approaches developed
by  Talagrand in~\cite{talagrand2006parisi} and earlier
Guerra in~\cite{guerra2003broken}. In particular, one creates
an interpolation between many instances $(H(\cdot,J_j), j\in [M])$ of the Hamiltonian arranged in a particular tree-like continuous  way. One then argues that a stable algorithm implemented on this ``continuous tree'' of instances results in a continuous (due to stability) tree of solutions,
the kind represented in Figure~\ref{fig:UTree-no-gap}. Then a separate argument is used to show that every ``continuous tree'' of 
spin configurations $(\sigma_j, j\in [M])$ will have at least one solution with energy at most $\eta_{p, \rm ALG}$, and thus (using an additional
concentration argument) stable algorithms cannot beat this threshold.  Putting it differently, $M$-tuples of spin configurations achieving 
simultaneously values
larger than $\eta_{p, \rm ALG}$ have associated ultrametric trees with gaps -- as represented by the red segments in Figure~\ref{fig:UTree-B-OGP}.

It is reasonable to believe that the scope of algorithms limited by the value $\eta_{p, \rm ALG}$ is much broader than the class
of stable algorithms, and quite possibly includes the class of all polynomial-time algorithms for these and similar problems.
 Such a result however is out of reach for now due to its very strong algorithmic complexity implications.

%\section{Quantum extensions. Perhaps discuss ML paper}

\bibliographystyle{foo}
%\bibliography{Bibliography/bibliography-07.2024}
%\bibliography{Bibliography/bibliography-07.2024}

%\bibliography{/Users/gamarnik/Documents/gamarnik/My_Papers/Bibliography/bibliography-12.2024}
\bibliography{bibliography-12.2024}

\begin{bibdiv}
\begin{biblist}

\bib{achlioptas2006solution}{inproceedings}{
      author={Achlioptas, Dimitris},
      author={Ricci-Tersenghi, Federico},
       title={On the solution-space geometry of random constraint satisfaction
  problems},
        date={2006},
   booktitle={Proceedings of the {T}hirty-{E}ighth {A}nnual {ACM} {S}ymposium
  on {T}heory of {C}omputing},
       pages={130\ndash 139},
}

\bib{bresler2022algorithmic}{inproceedings}{
      author={Bresler, Guy},
      author={Huang, Brice},
       title={The algorithmic phase transition of random {K-SAT} for low degree
  polynomials},
organization={{IEEE}},
        date={2022},
   booktitle={2021 {{IEEE}} 62nd {A}nnual {S}ymposium on {F}oundations of
  {C}omputer {S}cience ({{FOCS}})},
       pages={298\ndash 309},
}

\bib{chen2019suboptimality}{article}{
      author={Chen, Wei-Kuo},
      author={Gamarnik, David},
      author={Panchenko, Dmitry},
      author={Rahman, Mustazee},
       title={Suboptimality of local algorithms for a class of max-cut
  problems},
        date={2019},
     journal={The Annals of Probability},
      volume={47},
      number={3},
       pages={1587\ndash 1618},
}

\bib{coja2010better}{article}{
      author={Coja-Oghlan, Amin},
       title={A better algorithm for random {K-SAT}},
        date={2010},
     journal={SIAM Journal on Computing},
      volume={39},
      number={7},
       pages={2823\ndash 2864},
}

\bib{ding2022proof}{article}{
      author={Ding, Jian},
      author={Sly, Allan},
      author={Sun, Nike},
       title={Proof of the satisfiability conjecture for large $k$},
        date={2022},
     journal={Annals of Mathematics},
      volume={196},
      number={1},
       pages={1\ndash 388},
}

\bib{el2021optimization}{article}{
      author={El~Alaoui, Ahmed},
      author={Montanari, Andrea},
      author={Sellke, Mark},
       title={Optimization of mean-field spin glasses},
        date={2021},
     journal={The Annals of Probability},
      volume={49},
      number={6},
       pages={2922\ndash 2960},
}

\bib{FriezeIndependentSet}{article}{
      author={Frieze, A.},
       title={On the independence number of random graphs},
        date={1990},
     journal={Discrete Mathematics},
      volume={81},
       pages={171\ndash 175},
}

\bib{gamarnik2024hardness}{article}{
      author={Gamarnik, David},
      author={Jagannath, Aukosh},
      author={Wein, Alexander~S},
       title={Hardness of random optimization problems for boolean circuits,
  low-degree polynomials, and {L}angevin dynamics},
        date={2024},
     journal={SIAM Journal on Computing},
      volume={53},
      number={1},
       pages={1\ndash 46},
}

\bib{GamarnikKizildagPerkinsXu2022}{inproceedings}{
      author={Gamarnik, D.},
      author={K{\i}z{\i}lda{\u{g}}, E.},
      author={Perkins, W.},
      author={Xu, C.},
       title={Algorithms and barriers in the symmetric binary perceptron
  model.},
        date={2022},
   booktitle={{FOCS} 2022: 63rd {IEEE} symposium on foundations of computer
  science},
       pages={576\ndash 587},
}

\bib{gamarnik2022disordered}{article}{
      author={Gamarnik, David},
      author={Moore, Cristopher},
      author={Zdeborov{\'a}, Lenka},
       title={Disordered systems insights on computational hardness},
        date={2022},
     journal={Journal of Statistical Mechanics: Theory and Experiment},
      volume={2022},
      number={11},
       pages={114015},
}

\bib{gamarnik2014limits}{article}{
      author={Gamarnik, David},
      author={Sudan, Madhu},
       title={Limits of local algorithms over sparse random graphs},
        date={2017},
     journal={Annals of Probability},
      volume={45},
       pages={2353\ndash 2376},
}

\bib{guerra2003broken}{article}{
      author={Guerra, Francesco},
       title={Broken replica symmetry bounds in the mean field spin glass
  model},
        date={2003},
     journal={Communications in mathematical physics},
      volume={233},
       pages={1\ndash 12},
}

\bib{huang2022tight}{inproceedings}{
      author={Huang, Brice},
      author={Sellke, Mark},
       title={Tight {L}ipschitz hardness for optimizing mean field spin
  glasses},
organization={{IEEE}},
        date={2022},
   booktitle={2022 {IEEE} 63rd annual symposium on foundations of computer
  science ({FOCS})},
       pages={312\ndash 322},
}

\bib{karp1976probabilistic}{article}{
      author={Karp, Richard~M},
       title={The probabilistic analysis of some combinatorial search
  algorithms},
        date={1976},
     journal={Algorithms and complexity: New directions and recent results},
      volume={1},
       pages={1\ndash 19},
}

\bib{mezard2005clustering}{article}{
      author={M{\'e}zard, M.},
      author={Mora, T.},
      author={Zecchina, R.},
       title={Clustering of solutions in the random satisfiability problem},
        date={2005},
     journal={Physical Review Letters},
      volume={94},
      number={19},
       pages={197205},
}

\bib{montanari2021optimization}{article}{
      author={Montanari, Andrea},
       title={Optimization of the {S}herrington--{K}irkpatrick hamiltonian},
        date={2021},
     journal={SIAM Journal on Computing},
      number={0},
       pages={{FOCS}19\ndash 1},
}

\bib{parisi1980sequence}{article}{
      author={Parisi, Giorgio},
       title={A sequence of approximated solutions to the {SK} model for spin
  glasses},
        date={1980},
     journal={Journal of Physics A: Mathematical and General},
      volume={13},
      number={4},
       pages={L115},
}

\bib{rahman2017local}{article}{
      author={Rahman, Mustazee},
      author={Virag, Balint},
       title={Local algorithms for independent sets are half-optimal},
        date={2017},
     journal={The Annals of Probability},
      volume={45},
      number={3},
       pages={1543\ndash 1577},
}

\bib{subag2021following}{article}{
      author={Subag, Eliran},
       title={Following the ground states of full-{RSB} spherical spin
  glasses},
        date={2021},
     journal={Communications on Pure and Applied Mathematics},
      volume={74},
      number={5},
       pages={1021\ndash 1044},
}

\bib{talagrand2006parisi}{article}{
      author={Talagrand, Michel},
       title={The {P}arisi formula},
        date={2006},
     journal={Annals of mathematics},
       pages={221\ndash 263},
}

\end{biblist}
\end{bibdiv}

%/Users/gamarnik/Documents/gamarnik/My_Papers/Bibliography

\end{document}